\documentclass{amsart}
\numberwithin{equation}{section}
\usepackage[utf8]{inputenc}
\usepackage[T1]{fontenc}
\usepackage{lmodern}
\usepackage[margin=1.3in]{geometry}
\usepackage{setspace}
\usepackage[usenames, dvipsnames]{xcolor}

\usepackage{mathtools}
\usepackage{amssymb}
\usepackage{amsthm}

\usepackage[mathscr]{eucal}

\usepackage{tikz-cd}

\usepackage[backref=page, bookmarks=false]{hyperref}
\usepackage[capitalize, noabbrev]{cleveref}

\usepackage{microtype}

\newtheorem{thm}[equation]{Theorem}
\newtheorem*{thm*}{Theorem}
\newtheorem{lem}[equation]{Lemma}
\newtheorem{cor}[equation]{Corollary}
\newtheorem{prop}[equation]{Proposition}

\theoremstyle{definition}
\newtheorem{exm}[equation]{Example}
\newtheorem{defn}[equation]{Definition}

\newtheorem{ansatz}[equation]{Ansatz}

\theoremstyle{remark}
\newtheorem{rem}[equation]{Remark}

\crefname{thm}{Theorem}{Theorems}
\crefname{lem}{Lemma}{Lemmas}
\crefname{cor}{Corollary}{Corollaries}
\crefname{prop}{Proposition}{Propositions}
\crefname{ex}{Exercise}{Exercises}
\crefname{exm}{Example}{Examples}
\crefname{defn}{Definition}{Definitions}
\crefname{claim}{Claim}{Claims}
\crefname{rem}{Remark}{Remarks}
\crefname{fct}{Fact}{Facts}
\crefname{note}{Note}{Notes}
\crefname{ansatz}{Ansatz}{Ansatzes}

\hypersetup{
 colorlinks,
 linkcolor={red!50!black},
 citecolor={green!50!black},
 urlcolor={blue!80!black}
}

\newcommand{\term}{\emph}
\newcommand{\cc}{\mathrm{cc}}
\newcommand{\curv}{\mathrm{curv}}

\newcommand{\Z}{\mathbb Z}
\newcommand{\R}{\mathbb R}
\newcommand{\T}{\mathbb T}

\makeatletter
\def\instring#1#2{TT\fi\begingroup
  \edef\x{\endgroup\noexpand\in@{#1}{#2}}\x\ifin@}
\def\isuppercase#1{%
  \instring{#1}{ABCDEFGHIJKLMNOPQRSTUVWXYZ}%
}%
\newcommand{\C@lIfUpper}[1]{
 \if\isuppercase{#1}\mathscr{#1}%
 \else #1%
 \fi
}
\newcommand{\cat}[1]{\mathit{\@tfor\next:=#1\do{\C@lIfUpper{\next}}}}
\makeatother
\newcommand{\Man}{\cat{Man}}

\DeclarePairedDelimiter{\paren}{(}{)}
\DeclarePairedDelimiter{\abs}{\lvert}{\rvert}

\newcommand{\fg}{\mathfrak g}
\newcommand{\Sym}{\mathrm{Sym}}
\renewcommand{\Im}{\mathrm{Im}}
\renewcommand{\d}{\mathrm d}
\newcommand{\ud}{\,\d}

\newcommand{\Lhi}{L_{\mathrm{hi}}}
\newcommand{\Rhi}{R_{\mathrm{hi}}}
\newcommand{\Shp}{\Sh_{\mathrm{pu}}}
\newcommand{\Cyc}{\mathrm{Cyc}}
\newcommand{\Def}{\mathrm{Def}}
\newcommand{\id}{\mathrm{id}}
\newcommand{\pt}{\mathrm{pt}}

\newcommand{\Spc}{\cat{Sp}}
\newcommand{\fC}{\cat C}
\newcommand{\Sh}{\cat{Sh}}
\newcommand{\KU}{\mathit{KU}}
\newcommand{\KO}{\mathit{KO}}
\newcommand{\CW}{\mathit{CW}}
\newcommand{\SO}{\mathrm{SO}}
\newcommand{\U}{\mathrm U}
\renewcommand{\O}{\mathrm{O}}
\newcommand{\ch}{\mathrm{ch}}
\newcommand{\ph}{\mathrm{ph}}
\newcommand{\Zsm}{Z^{\mathrm{sm}}}
\newcommand{\Csm}{C^{\mathrm{sm}}}

\newcommand{\CS}{\mathrm{CS}}

\newcommand{\many}[1]{#1\dotsb#1}

\title{Differential cohomology}
\author{Arun Debray}
\address{Department of Mathematics, Purdue University, 150 N. University Street, West Lafayette, IN 47907-2067}
\email{adebray@purdue.edu}
\date{\today}

\keywords{Differential cohomology, Chern-Weil theory, Cheeger-Simons differential characters, Deligne cohomology}

\begin{document}
\maketitle

\begin{abstract}
We give an overview of differential cohomology from the point of view of algebraic topology. This includes a survey
of several different definitions of differential cohomology groups, a discussion of differential characteristic
classes, an introduction to differential generalized cohomology theory, and some applications in physics.
\end{abstract}

\tableofcontents

\setcounter{section}{-1}
\section{Introduction}
	It is a truth universally acknowledged that a closed differential form, in possession of an interpretation as a
gauge field in a quantum field theory, must be in want of an integral refinement. This refinement manifests the
quantum nature of the physical theory: that quantities in the theory are ``quantized,'' meaning that in some system
of units they are integers, not arbitrary real numbers. The mathematical incarnation of this theory of closed forms
with integrality data is called \term{differential cohomology}; the objective of this article is to survey this
theory, including several different approaches to the basic definitions, some useful constructions in the theory,
and some applications.

The basic data is as follows. For $M$ a smooth manifold, there are differential cohomology groups $\check
H^k(M;\Z)$  equipped with a \term{characteristic class map} $\cc\colon\check H^k(M;\Z)\to H^k(M;\Z)$ and a
\term{curvature map} $\curv\colon\check H^k(M;\Z)\to\Omega^k(M)_{\mathit{c\ell}}$; there is a sense in which
$\check H^k(M;\Z)$ is the universal object classifying data of an integral cohomology class (its characteristic
class), a closed differential form (its curvature), and an identification of the two induced de Rham cohomology
classes. The first construction of $\check H^k(M;\Z)$ was given by Cheeger-Simons~\cite{CS85}, and since then
many constructions, concrete and abstract, have appeared; we will survey several in \S\ref{defns}.

It is a general rule of thumb that ordinary cohomology is to topological objects as differential cohomology is to
geometric ones. For example, $H^2(M;\Z)$ classifies complex line bundles $L\to M$, and $\check H^2(M;\Z)$
classifies complex line bundles with connection. The characteristic class and curvature maps capture the first
Chern class of the line bundle, resp.\ the curvature of the connection.

As differential cohomology feels like ordinary cohomology, but upgraded, one can ask which facts about ordinary
cohomology upgrade to differential cohomology. The answer is that quite a lot of them do, including integration
along the fiber of a relatively oriented bundle of smooth manifolds. In addition, large parts of the theory of
characteristic classes lift to differential cohomology, and even enhance: the differential cohomology refinement of
the Chern-Weil map contains the information of Chern-Simons invariants, for example. In \S\ref{char_classes} we discuss
this and other aspects of characteristic classes in differential cohomology. Another avenue for analogy with
ordinary cohomology is the prospect of differential refinements of generalized cohomology theories.  These too
exist, and we will discuss theory and examples in \S\ref{gencoh}.

In \S\ref{s_physics}, we discuss applications of differential cohomology in theoretical physics: quantization of
abelian gauge fields.

Finally, in \S\ref{further_reading}, we give some suggestions for further reading.

\section{Definitions}
	\label{defns}
	Before the main content of this section, where we survey several different models for differential cohomology, let
us begin with some basic key facts about these groups.

Differential cohomology is a theory assigning to each smooth manifold $M$ a series of abelian Fréchet Lie groups
$\check H^k(M;\Z)$ which are \emph{not} homotopy invariants of $M$.
\begin{itemize}
	\item $\check H^1(M;\Z)$ is naturally isomorphic to the set of functions $M\to \R/\Z$.
	\item $\check H^2(M;\Z)$ is naturally isomorphic to the isomorphism classes of complex line bundles on $M$ with
	connection.
	\item $\check H^*(M;\Z)$ comes with a cup product making it into a ring.
	\item There is a map $\cc\colon\check H^k(M;\Z)\to H^k(M;\Z)$, called the \term{characteristic class map}. On
	$\check H^2$, this sends a line bundle with connection to its first Chern class.
	\item There is a map $\curv\colon \check H^k(M;\Z)\to\Omega^k(M)_{\mathit{c\ell}}$ (i.e.\ to closed $k$-forms),
	called the \term{curvature map}. On $\check H^2$, this sends a line bundle with connection to its curvature.
\end{itemize}

\subsection{Cheeger-Simons' differential characters}
Differential characters are for the reader who sees de Rham cohomology's philosophy of form over function and
thinks, ``why can't I have both?'' They were the first definition of differential cohomology to appear, and have
the feel of singular cohomology.
\begin{defn}[{Cheeger-Simons~\cite[\S 1]{CS85}}]
Let $M$ be a smooth manifold and write $\Csm_k(M)$, resp.\ $\Zsm_k(M)$, for the abelian groups of smooth
$k$-chains, resp.\ $k$-cycles on $M$. A \term{differential character} of degree $n$ on $M$ is a homomorphism
$\chi\colon \Zsm_{n-1}(M)\to\R/\Z$ such that there exists $\omega\in\Omega^n(M)$ such that for all $C\in \Csm_{n-1}(M)$,
\begin{equation}
	\chi(\partial c) = \int_C \omega(\chi)\bmod \Z.
\end{equation}
The degree-$n$ \term{differential cohomology} of $M$, denoted $\check H^n(M;\Z)$, is the group of degree-$n$
differential characters.
\end{defn}
There is a unique $\omega$ satisfying this definition for a given $\chi$, and $\omega$ is always a closed form. The
curvature map sends $\chi\mapsto\omega$.

The characteristic class map has a slightly more elaborate definition. Since $\Zsm_{n-1}(M)$ is a free abelian
group and $\R\to\R/\Z$ is an epimorphism, $\chi\colon\Zsm_{n-1}(M)\to\R/\Z$ lifts to a homomorphism
$\widetilde\chi\colon\Zsm_{n-1}(M)\to\R$. Now define $I(\widetilde\chi)\colon\Csm_{n-1}(M)\to\Z$ by
\begin{equation}
	C\longmapsto -\widetilde\chi(\partial C) + \int_C \curv(\chi).
\end{equation}
One can show this is indeed $\Z$-valued, and since $\curv(\chi)$ is a closed form, this is a cocycle. The
characteristic class morphism sends $\chi\mapsto [I(\widetilde\chi)]$, which can be shown to not depend on the
choice of lift $\widetilde\chi$.
\begin{rem}
Our indexing convention differs from Cheeger-Simons' original convention; we follow the standard convention in the
field of differential cohomology, so that the characteristic class and curvature morphisms preserve the degree.
\end{rem}

\subsection{Deligne cohomology}
Deligne cohomology refers to a sheaf cohomology model for differential cohomology. Deligne first studied this model
in an algebro-geometric setting in~\cite{Del71}. Brylinski~\cite{Bry93} was the first to consider this model on
smooth manifolds.

Throughout this article, we make the category $\Man$ of manifolds and smooth functions into a site in which
the coverings are surjective submersions, and we define a few sheaves on this site.
\begin{itemize}
	\item Given a Lie group $A$, we let $\underline A$ denote the sheaf of smooth $A$-valued functions; $A$ without
	underline denotes the sheaf of locally constant $A$-valued functions, i.e.\ smooth functions for the discrete
	Lie group structure on $A$.
	\item The \term{sheaf of differential $k$-forms}, denoted $\Omega^k$, sends a manifold $M$ to the real vector
	space $\Omega^k(M)$ of $k$-forms on $M$. Thus we have an isomorphism $\upsilon\colon \underline
	\R\overset\cong\to\Omega^0$.
\end{itemize}
We will also take chain complexes of sheaves on $\Man$. The categories of chain complexes of sheaves on $\Man$ and
sheaves of chain complexes on $\Man$ are isomorphic; given a sheaf of chain complexes $\mathcal F^\bullet$ on
$\Man$ and a smooth manifold $M$, $H^*(M;\mathcal F^\bullet)$ refers to the \term{hypercohomology} of $M$ valued in
$C$; that is, form the double complex $C^p(M; \mathcal F^q)$ and take the cohomology with respect to the total
differential.\footnote{For an abelian group $A$, $H^*(M;A)$ is a priori ambiguous --- do we mean singular
$A$-cohomology of $M$ or the sheaf cohomology of $M$ with respect to the sheaf $A$? Fortunately, these two
cohomology theories are naturally isomorphic.}
\begin{defn}[{Deligne~\cite[\S 2.2]{Del71}}]
The \term{Deligne complex} $\Z(n)$ is the chain complex of sheaves
\begin{subequations}
\begin{equation}
	\Z(n) \coloneqq \paren{0\longrightarrow \Z\longrightarrow \Omega^0 \many\longrightarrow
	\Omega^{n-1}\longrightarrow 0}.
\end{equation}
Here the map $\Z\to\Omega^0$ is the inclusion of $\Z$-valued functions into $\R$-valued functions combined with the
isomorphism $\upsilon \colon \underline\R\overset\cong \to\Omega^0$.
\begin{rem}
Etymologically, the Deligne complex $\Z(n)$ is related to the ``Tate twist'' that is also often denoted $\Z(n)$,
but the two are not equivalent. For this reason, some authors denote the Deligne complex something like
$\Z(n)_{\mathcal D}$.
\end{rem}
One also sees the complexes
\begin{align}
	\R(n) &\coloneqq \paren{0\longrightarrow \R\overset{\widetilde\upsilon}{\longrightarrow} \Omega^0
	\many\longrightarrow \Omega^{n-1}\longrightarrow 0},
\intertext{where $\widetilde\upsilon$ is the inclusion of locally constant functions into all functions followed by
$\upsilon$, and}
	\T(n) &\coloneqq \paren{0\longrightarrow \underline \T \overset\varphi\longrightarrow \Omega^1
	\many\longrightarrow \Omega^n\longrightarrow 0},
\end{align}
where $\varphi\coloneqq (1/2\pi i)\d\log$ and $\d\log\colon\T\to i\Omega^1$ is the morphism sending a $\T$-valued
function $f$ to the form $(1/f)\ud f$.
\end{subequations}
\end{defn}
\begin{prop}[{Brylinski~\cite[Proposition 1.5.7]{Bry93}}]
For any manifold $M$, there is a natural isomorphism $H^n(M;\Z(n))\cong\check H^n(M)$.
\end{prop}
Thus the ``diagonally graded'' Deligne cohomology groups are differential cohomology groups. The ``off-diagonal''
groups $H^k(M;\Z(n))$ for $k\ne n$ are isomorphic to singular cohomology valued in $\Z$ (if $k>n$) or $\R/\Z$ (if
$k<n$) (see~\cite[Theorem 1.5.3]{Bry93} and~\cite[\S 3.2]{HS05}), so appear uninteresting at first glance, but they
attain interesting values on certain stacks; see \S\ref{off_diagonal}.

$\R(n)$ and $\T(n)$ are also familiar: $\R(n)$ is isomorphic to the sheaf of closed $n$-forms considered as a
complex concentrated in degree $n$, and $\T(n)\simeq \Z(n+1)[-1]$~\cite[Remark 3.6]{BM94}; the proof of the latter
essentially amounts to the weak equivalence of the complexes $0\to\Z\to\underline\R\to 0$ and $0\to
0\to\underline\T\to 0$.

In this model for differential cohomology, $\curv\colon H^n(M;\Z(n))\to \Omega_{\mathit{c\ell}}^n(M)$ is the map
$\Z(n)\to\R(n)$ induced by the inclusion $\Z\hookrightarrow\R$, together with the identification of
$\R(n)$-cohomology with closed $n$-forms. The characteristic class map is the effect on cohomology of the
\term{truncation map} $t\colon \Z(n)\to\Z$ defined by quotienting $\Z(n)$ by the subcomplex of sheaves in positive
homological degrees.

Harvey-Lawson's model for differential cohomology in terms of ``sparks''~\cite{HL06} has a similar feel to Deligne
cohomology, though with a cocycle model.
\subsection{Hopkins-Singer's homotopy pullback model}
Hopkins-Singer's approach to differential cohomology~\cite[\S 3]{HS05} begins with the following observation.
\begin{lem}[{Hopkins-Singer~\cite[\S 3.2]{HS05}, Bunke-Nikolaus-Völkl~\cite[\S 4.1]{BNV16}}]\label{HP_lem}
The truncation maps $t\colon\Z(n)\to\Z$ and $\R(n)\to\R$ participate in a homotopy pullback square
\begin{equation}
\label{homotopy_pullback_square}
\begin{tikzcd}
	{\Z(n)} & \Z \\
	{\R(n)} & \R,
	\arrow["t", from=2-1, to=2-2]
	\arrow["t", from=1-1, to=1-2]
	\arrow[from=1-1, to=2-1]
	\arrow[from=1-2, to=2-2]
		\arrow["\lrcorner"{anchor=center, pos=0.125}, draw=none, from=1-1, to=2-2]
\end{tikzcd}\end{equation}
where the vertical arrows are induced by the usual inclusion $\Z\to\R$.
\end{lem}
This expresses the idea that differential cohomology is a homotopy pullback of closed differential forms and
integral cohomology. Hopkins-Singer provide an explicit cocycle model for this homotopy pullback.
\begin{defn}[{Hopkins-Singer~\cite[\S 3.2]{HS05}}]
Let $\check C(q)^\bullet(M)$ be the cochain complex given by
\begin{equation}
	\check C(q)^n(M) \coloneqq\begin{cases}
		C^n(M;\Z)\times C^{n-1}(M;\R) \times \Omega^n(M), &n\ge q\\
		C^n(M;\Z)\times C^{n-1}(M;\R), &n<q,
	\end{cases}
\end{equation}
with differential given by, when $n\ge q$,
\begin{subequations}
\begin{equation}
	d(c, h, \omega) \coloneqq (\delta c, \omega - c - \delta h, \d\omega)
\end{equation}
and when $n < q$,
\begin{equation}
	d(c, h) \coloneqq\begin{cases}
		(\delta c, -c-\delta h, 0), &n = q-1\\
		(\delta c, -c-\delta h), &n < q-1.
	\end{cases}
\end{equation}
\end{subequations}
The degree-$n$ \term{differential cohomology} of $M$ is $\check H^n(M;\Z)\coloneqq H^n(\check C(n)^\bullet(M))$.
\end{defn}
If $(c,h,\omega)$ is an $n$-cocycle, then $c$ and $\omega$ are both closed; the characteristic class map sends
$(c,h,\omega)$ to the class of $c$, and the curvature map sends $(c,h,\omega)\mapsto\omega$.
\subsection{Simons-Sullivan's hexagon}
Simons-Sullivan~\cite{SS08} produced a property of differential cohomology that uniquely characterizes it, in terms
of a hexagon-shaped diagram.
\begin{thm}[{Simons-Sullivan~\cite[Theorem 1.1]{SS08}}]
\label{hexagon_thm}
Let $\widehat H^*$ be a functor from manifolds to graded abelian groups, and suppose $\widehat H$ is equipped with natural
transformations
\begin{enumerate}
	\item $i_1\colon H^{n-1}(\text{--};\R/\Z)\to \widehat H^n$,
	\item $i_2\colon \Omega^{n-1}(\text{--})/\Im(\d)\to\widehat H^k$,
	\item $\delta_1\colon\widehat H^n\to\Omega^n_{\mathit{c\ell}}$, and
	\item $\delta_2\colon \widehat H^n\to H^n(\text{--};\Z)$
\end{enumerate}
such that the following diagram commutes:
\begin{equation}\label{the_hexagon}
\begin{tikzcd}
	& {H^{n-1}(\text{--};\R/\Z)} && {H^n(\text{--};\Z)} \\
	{H^{n-1}(\text{--};\R)} && {\widehat H^n} && {H^n(\text{--};\R)}, \\
	& {\Omega^{n-1}(\text{--})/\Im(\d)} && {\Omega^n_{\mathit{c\ell}}}
	\arrow[from=2-1, to=1-2]
	\arrow[from=2-1, to=3-2]
	\arrow[from=1-2, to=1-4]
	\arrow[from=1-4, to=2-5]
	\arrow["\d", from=3-2, to=3-4]
	\arrow[from=3-4, to=2-5]
	\arrow["{i_1}"{description}, from=1-2, to=2-3]
	\arrow["{i_2}"{description}, from=3-2, to=2-3]
	\arrow["{\delta_1}"{description}, from=2-3, to=3-4]
	\arrow["{\delta_2}"{description}, from=2-3, to=1-4]
\end{tikzcd}\end{equation}
where the topmost arrows are the Bockstein long exact sequence associated to $0\to\Z\to\R\to\R/\Z\to 0$, and the
bottommost arrows come from the de Rham theorem. Then $\widehat H^*\cong\check H^*$, $\delta_1$ and $\delta_2$ are
respectively the curvature and characteristic class maps, and $i_1$ and $i_2$ are their kernels.
\end{thm}
So this diagram, the \term{differential cohomology hexagon}, contains quite a bit of information: the topmost
arrows are a long exact sequence, the bottommost arrows are another long exact sequence, and the diagonals extend
to a short exact sequence. Moreover, both squares, when lifted to the level of sheaves of complexes on $\Man$, are
homotopy pullback squares.

See Stimpson~\cite{Sti11} for another axiomatic characterization of differential cohomology.

\subsection{A few basic properties of differential cohomology}
We conclude this section with a few elementary properties of differential cohomology.
\subsubsection{Higher gerbes with connection}
Recall that $\check H^1(M;\Z)$ is the group of functions $M\to \R/\Z$ -- or equivalently, $M\to\T$, and that
$\check H^2(M;\Z)$ is the isomorphism classes of complex line bundles with connection. These can be thought of as
categorifications of $\T$-valued functions, suggesting that higher differential cohomology groups ought to
represent categorifications of the notion of line bundle with connection. This is correct: these higher-categorical
objects are called \term{gerbes with connection}. See Brylinski~\cite{Bry93} and the references therein for more
information.
\subsubsection{Cup product}
The differential cohomology groups $\check H^*(M;\Z)$ jointly carry a ring structure, which the characteristic
class map sends to the ordinary cup product and the curvature map sends to wedge product. There are various
different ways to construct this cup product dating back to Cheeger-Simons' original construction~\cite{CS85} of
differential cohomology.

\subsubsection{Integration along the fiber}
Suppose $E\to B$ is a submersion of smooth manifolds with fiber $F$ and an oriented vertical tangent bundle ---
i.e.\ exactly the conditions needed to have integration along the fiber in ordinary cohomology. Then, there is also
a differential-cohomology integration along the fiber:
\begin{equation}
	\int_F\colon H^k(E; \Z(n)) \longrightarrow H^{k-\dim(F)}(B; \Z(n-\dim(F))).
\end{equation}
See Hopkins-Singer~\cite[\S 2.4]{HS05} for a construction of this map.
\subsubsection{Cohomology operations}
Grady-Sati~\cite{GS18b, GS18a} have lifted primary and secondary cohomology operations to differential cohomology.

\section{Differential characteristic classes}
	\label{char_classes}
	Characteristic classes are a place where differential cohomology shines: the analogy with ordinary cohomology is
close enough to help both intuition and proofs, yet differential characteristic classes engender new phenomena,
including geometric invariants such as Chern-Simons invariants. Chern-Weil theory is the key to the story, so we
begin with that, and then explain how it manifests in differential cohomology.

\subsection{Chern-Weil theory}
In this subsection only, we will let $\Omega_M^*(V)$ denote differential forms on the manifold $M$ valued in the
vector space $V$, i.e.\ sections of the bundle $\Lambda^*(T^*M\otimes V)$. For example, a connection $\Theta$ on a
principal $G$-bundle $P\to M$ is a form in $\Omega_P^1(\fg)$, where $\fg$ is the Lie algebra of $G$, and the
curvature of $\Theta$ is an element of $\Omega_P^2(\fg)$.

Consider the algebra $\Sym^*(\fg^\vee)^G$, i.e.\ the $G$-invariants of the algebra of polynomial functions
$\fg\to\R$. Here the $G$-action is induced from the adjoint action of $G$ on $\fg$. Given $f\in\Sym^k(\fg^\vee)^G$,
so that $f$ is a degree-$k$ polynomial, together with a principal $G$-bundle with connection, we will build a
closed $2k$-form whose de Rham class is a characteristic class for $G$-bundles.

Let $\pi\colon P\to M$ be a principal $G$-bundle with connection $\Theta$ and curvature $\Omega\in\Omega_P^2(\fg)$.
Wedge together $k$ copies of $\Omega$ to produce $\Omega^{\wedge k}\in \Omega_P^{2k}(\fg^{\otimes k})$, then apply
the polynomial $f$ to obtain $f(\Omega^{\wedge k})\in\Omega_P^{2k}(\R)$. Because $f$ is $\mathrm{Ad}$-invariant,
this form descends to a form $\CW(P, \Theta, f)\in\Omega^{2k}_M(\R)$, called the \term{Chern-Weil form} of $P$,
$\Theta$, and $f$. The key properties of Chern-Weil forms are:
\begin{enumerate}
	\item $\CW(P, \Theta, f)$ is a closed form,
	\item the de Rham class of $\CW(P, \Theta, f)$ depends on $P$ but not on the choice of $\Theta$, and
	\item holding $f$ fixed, $\CW(P,\Theta,f)$ is natural in $(P,\Theta)$.
\end{enumerate}
In fact, the Chern-Weil construction defines an isomorphism $\Sym^*(\fg^\vee)^G\to H^*(BG;\R)$ for any compact Lie
group $G$.

\subsection{Differential cohomology lifts of Chern-Weil forms}
Suppose that the de Rham class of a Chern-Weil form $\CW(P,\Theta, f)$ is in the image of the map
$H^*(\text{--};\Z)\to H^*(\text{--};\R)$. Then we have, at least at a heuristic level, the data of a differential
cohomology class: a closed form and an integral cohomology class with identified values in de Rham cohomology. Is
there a lift to differential cohomology? Cheeger-Simons~\cite{CS85} showed the answer is yes, and
Bunke-Nikolaus-Völkl~\cite{BNV16} showed the naturality properties of these classes allow one to work universally
with the classifying stack $B_\nabla G$ of principal $G$-bundles with connection.\footnote{By a \term{stack} we
mean a simplicial sheaf on $\Man$. See Freed-Hopkins~\cite{FH13} for more information.}
\begin{thm}[{Cheeger-Simons~\cite[Theorem 2.2]{CS85}, Bunke-Nikolaus-Völkl~\cite[\S 5.2]{BNV16}}]
\label{diffCW}
Let $G$ be a compact Lie group and $c\in H^{2k}(BG;\Z)$. Let $f\in\Sym^k(\fg^\vee)^G$ be the polynomial uniquely
characterized by asking for the de Rham class of $\CW(P,\Theta, f)$ to equal $c(P)$ in $H^{2k}(BG;\R)$. Then, there
is a unique class $\check c\in \check H^{2k}(B_\nabla G;\Z)$ whose characteristic class is $c$ and whose curvature
form is the Chern-Weil form.
\end{thm}
That is, Chern-Weil theory produces characteristic classes in differential cohomology depending on a principal
bundle and a connection.
\begin{exm}
For $G = \O_n$, $\SO_n$, or $\U_n$, we obtain characteristic classes of vector bundles with certain classes of
connections.
\begin{enumerate}
	\item The Pontrjagin classes $p_i\in H^{4i}(BO_n;\Z)$ of a real vector bundle lift to \term{differential
	Pontrjagin classes} $\check p_i\in \check H^{4i}(B_\nabla \O_n;\Z)$ of vector bundles equipped with a metric and
	a compatible connection. See Brylinski-McLaughlin \cite{BM96} and Grady-Sati \cite[Proposition 3.6]{GS21} for
	additional constructions of these classes.
	\item The Chern classes $c_i\in H^{2i}(B\U_n;\Z)$ of a complex vector bundle lift to \term{differential Chern
	classes} $\check c_i\in\check H^{2i}(B_\nabla \U_n;\Z)$ of complex vector bundles equipped with a Hermitian
	metric and a compatible connection.
	See Several authors construct differential Chern classes by other methods, including Brylinski-McLaughlin
	\cite{BM96}, Berthomieu \cite{Ber10}, Bunke \cite{Bun10,Bunke:notes}, and Ho \cite{Ho15} for additional
	constructions of these classes.
	\item The Euler class $e\in H^{2k}(B\SO_{2k};\Z)$ of an oriented real rank-$2k$ vector bundle lifts to a
	\term{differential Euler class} $\check e\in \check H^{2k}(B_\nabla \SO_{2k};\Z)$ of such vector bundles
	equipped with a metric and a compatible connection. See Brylinski-McLaughlin \cite{BM96} and Bunke
	\cite[Example 3.85]{Bunke:notes} for additional constructions of $\check e$.
\end{enumerate}
For differential Pontrjagin and Chern classes but \emph{not} the differential Euler class, one can relax the
condition of compatibility with the metric. See~\cite[Remark 14.1.13]{ADH21}. These classes also satisfy a Whitney
sum formula, as described in~\cite[\S 14.2]{ADH21}.
\end{exm}
Fiorenza-Schreiber-Stasheff~\cite{FSS12} generalize this story to higher groups.
\subsection{Chern-Simons invariants}
Choose a class $c\in H^{2k}(BG;\Z)$ and let $\check c\in\check H^{2k}(B_\nabla G;\Z)$ be its differential
refinement as guaranteed by \cref{diffCW}. If $M$ is a closed, oriented $2k$-manifold together with a principal
$G$-bundle $P\to M$ with connection $\Theta$, the data $(P,\Theta)$ pull $\check c$ back to a class $\check
c(P,\Theta)\in\check H^{2k}(M;\Z)$. Since $M$ is oriented, we can integrate:
\begin{equation}
	\int_M \check c(P,\Theta) \in\check H^0(\pt;\Z)\cong\Z.
\end{equation}
This integral is not so interesting: we just recover $\int_M c(P)$, as if we had never entered the world of
differential cohomology. But there is something better we can do: since $\check H^1(\pt;\Z)\cong\R/\Z$, we can
integrate on a $(2k-1)$-manifold $N$ with principal $G$-bundle $P$ and connection $\Theta$:
\begin{equation}
	\int_N\check c(P,\Theta)\in\check H^1(\pt;\Z)\cong\R/\Z.
\end{equation}
This number turns out to be more interesting --- it recovers the Chern-Simons invariant.
\begin{defn}[{Chern-Simons~\cite{CS74}}]
Choose $f\in\Sym^k(\fg^\vee)^G$ and let $P\to M$ be a principal $G$-bundle with two connections $\Theta_0$ and
$\Theta_1$ on it. Since the space of connections is convex, we can let $\Theta_t\coloneqq (1-t)\Theta_0 +
t\Theta_1$ for $t\in[0,1]$; these connections stitch together to a connection $\overline\Theta$ on $[0,1]\times M$,
with curvature $\overline\Omega$.

The \term{Chern-Simons form} of $P$, $\Theta_1$, and $\Theta_2$ is
\begin{equation}
	\CS_f(\Theta_1, \Theta_2) = \int_0^1 f(\overline\Omega)\in \Omega^{2k-1}(M).
\end{equation}
Given a bundle $p\colon P\to M$ with one connection $\Theta$, the Chern-Simons form $\CS_f(\Theta)$ is defined by
pulling $(P,\Theta)$ back along $P\to M$, then computing $\CS_f(\Theta^{\mathrm{triv}}, p^*\Theta)$, where
$\Theta^{\mathrm{triv}}$ is the connection coming from the tautological trivialization of $p^*P\to P$.
\end{defn}
\begin{prop}
\label{CS_prop}
With $c$ and $\check c$ as above, let $f$ be the invariant polynomial whose Chern-Weil form is $\curv(\check c)$,
and let $i_2\colon \Omega^{2k-1}(M)/\Im(\d)\to\check H^{2k}(M;\Z)$ be as in~\eqref{the_hexagon} (i.e.\ the kernel of
the characteristic class map). Then, for any principal $G$-bundle $\pi\colon P\to M$ with connection $\Theta$,
\begin{equation}
	i_2(\CS_f(\Theta)) = \pi^*\check c(P,\Theta)\in \check H^{2k}(P; \Z).
\end{equation}
\end{prop}
This was known to Chern-Simons~\cite{CS74} albeit not stated explicitly there; see~\cite[Proposition 19.1.9]{ADH21}
for a proof.

To more explicitly relate \cref{CS_prop} to the integration story we began with, fix $2k = 4$ and $G$ to be a
simple, simply connected Lie group, so that $BG$ is $3$-connected and any principal $G$-bundle $P\to N$ over a
$3$-manifold $N$ admits a section. Choose a section, and call it $s\colon N\to P$. Then
\begin{equation}
	\int_N s^*\CS_f(\Theta) = \int_N \check c(P,\Theta)\in\R/\Z.
\end{equation}
The left-hand side is a priori $\R$-valued, but depends on the section; the value in $\R/\Z$ is independent of the
choice of $s$.
\begin{rem}
For $\check c_1$, $\exp\left(2\pi i\int_{S^1}\check c_1(P,\Theta)\right)$ computes the holonomy of the connection
$\Theta$ around $S^1$.
\end{rem}

\subsection{Off-diagonal characteristic classes}
\label{off_diagonal}
Let $G$ be a Lie group and $B_\bullet G$ be the classifying stack of principal $G$-bundles.
Paralleling work of Beĭlinson~\cite[\S 1.7]{Bei84}, Bloch~\cite{Blo78}, Soulé~\cite{Sou89}, Brylinski~\cite{Bry99,
Bry99a}, and Dupont-Hain-Zucker~\cite{DHZ00} in algebraic geometry, people have lifted characteristic classes of
principal $G$-bundles to the ``off-diagonal'' Deligne cohomology groups $H^{2q}(B_\bullet G; \Z(q))$, beginning
with work of Bott~\cite{Bot73} calculating $H^*(B_\bullet G;\Omega^q)$ and of Shulman~\cite{Shu72} and
Bott-Shulman-Stasheff~\cite{BSS76} on $H^*(B_\bullet G; \Omega^{\ge q})$; these calculations were interpreted in
differential cohomology by Waldorf~\cite{Wal10} and in~\cite[Chapters 15--17]{ADH21}. One key result is a lift of
the Chern-Weil map.
\begin{thm}[Bott~\cite{Bot73}, Hopkins]
\label{big_off_diagonal_thm}
Let $G$ be a Lie group with $\pi_0(G)$ finite, let $i\colon K\hookrightarrow G$ be the inclusion of the maximal
compact subgroup of $G$, and let $\fg$ and $\mathfrak k$ be the Lie algebras of $G$ and $K$, respectively. Then
there is a commutative diagram
\begin{equation}\begin{tikzcd}
	{H^{2n}(B_\bullet G;\Z(n))} & {H^{2n}(BG;\Z)} & {H^{2n}(BK;\Z)} \\
	{\Sym^n(\fg^\vee)^G} & {\Sym^n(\mathfrak k^\vee)^K} & {H^{2n}(BK;\R)}
	\arrow["t", from=1-1, to=1-2]
	\arrow["{i^*}", "\cong"', from=1-2, to=1-3]
	\arrow[from=1-3, to=2-3]
	\arrow["{\mathit{CW}}", "\cong"', from=2-2, to=2-3]
	\arrow["{i^*}", from=2-1, to=2-2]
	\arrow[from=1-1, to=2-1]
	\arrow[from=1-2, to=2-2]
	\arrow["\lrcorner"{anchor=center, pos=0.125}, draw=none, from=1-1, to=2-2]
\end{tikzcd}\end{equation}
where the left-hand square is a pullback square and $\mathit{CW}$ is the usual Chern-Weil isomorphism.
\end{thm}
Bott~\cite{Bot73} proved a related result; this reinterpretation is due to Hopkins, and a proof can be found
in~\cite[Corollaries 16.2.4, 16.2.5]{ADH21}.

Hence if $G$ is compact, the truncation map $t\colon H^{2n}(B_\bullet G;\Z(n))\to H^{2n}(BG;\Z)$ is an isomorphism.
For noncompact $G$, \cref{big_off_diagonal_thm} allows one to use information on $\Sym^*(\fg^\vee)^G$ to gain
leverage on characteristic classes in $H^{2n}(B_\bullet G;\Z(n))$; see~\cite[Chapter 17]{ADH21} and~\cite[\S
3]{DLW23} for examples of this technique.

If $H$ is a Fréchet Lie group, there is a natural isomorphism from $H^3(B_\bullet H;\Z(1))$ to the abelian group of
Fréchet Lie group central extensions of $H$ by $\T$~\cite[Corollary 18.3.2]{ADH21}. Thus one can construct such extensions by
using \cref{big_off_diagonal_thm} to construct an off-diagonal characteristic class, then move it into
$H^3(\text{--};\Z(1))$ using some sort of transgression map. Work of Brylinski-McLaughlin~\cite[\S 5]{BM94} shows
how to use this to construct the Kac-Moody central extensions of loop groups of compact simple Lie
groups,\footnote{Brylinski-McLaughlin did not have \cref{big_off_diagonal_thm} available, so constructed their
off-diagonal differential characteristic classes a different way, using objects called \term{multiplicative bundle
gerbes}. \Cref{big_off_diagonal_thm} gives an alternative to that part of their proof.} and~\cite{DLW23} uses this
approach to construct the Virasoro central extensions of $\mathrm{Diff}^+(S^1)$.

\section{Differential generalized cohomology}
	\label{gencoh}
	Generalized cohomology theories such as $K$-theory and cobordism have long been an important ingredient in the
algebraic topologist's toolbox. In differential cohomology, analogous theories were motivated by ideas in string
theory, before more recent work studying all such ``differential generalized cohomology theories''\footnote{One
hears both ``differential generalized cohomology theory'' and ``generalized differential cohomology theory.'' In
this article, we favor the former: the way these theories have been studied in the literature generally treats them
as differential analogues of generalized cohomology theories, rather than generalizations of ordinary differential
cohomology. For example, one does not often see Eilenberg-Steenrod-type axioms for differential generalized
cohomology theories.} from a homotopical point of view. In this section, we will begin with the general theory
in~\S\ref{gen_theory}, then turn to examples in~\S\ref{gen_examples}.

Differential generalized cohomology theories were first proposed by Freed~\cite[\S 1]{Fre00}, who sketched a
definition. Hopkins-Singer~\cite[\S 4]{HS05} provided the first comprehensive treatment of differential generalized
cohomology. Bunke-Nikolaus-Völkl~\cite{BNV16} and Schreiber~\cite{Sch13} provide additional, more homotopical
treatments; in this section, we will follow Bunke-Nikolaus-Völkl's account.
\subsection{Differential generalized cohomology theories and sheaves on manifolds}
\label{gen_theory}
Let $\Spc$ denote the $\infty$-category of spectra, and for any presentable $\infty$-category $\fC$, such as
$\Spc$, let $\Sh(\Man, \fC)$ denote the $\infty$-category of $\fC$-valued sheaves on $\Man$. These are the functors $\mathcal F\colon\Man^{\cat{op}}\to\fC$ whose restriction to each
manifold is a sheaf in the usual sense.
\begin{defn}[Bunke-Nikolaus-Völkl~\cite{BNV16}]
A \term{differential generalized cohomology theory} is a cohomology theory on $\Man$ given by the sheaf cohomology
of some sheaf in $\Sh(\Man, \Spc)$.
\end{defn}
That is, generalized cohomology theories are to $\Spc$ as differential generalized cohomology theories are to
$\Sh(\Man, \Spc)$. Much of the theory in this section works with target an arbitrary presentable $\infty$-category
$\fC$ in place of $\Spc$; see~\cite{Sch13, ADH21} for more information.

Generalized differential cohomology theories are in general not homotopy-invariant. One easy example is
$H\Omega^k$, given by composing the sheaf of differential $k$-forms with the Eilenberg-Mac Lane functor.
$H\Omega^k$-cohomology is nontrivial on $\R^k$.
\begin{defn}
A sheaf $\mathcal F\in\Sh(\Man, \Spc)$ is \term{homotopy invariant}, or \term{concordance-invariant}, or
\term{$\R$-invariant}, if for every map of manifolds $f\colon M\to N$
that is a homotopy equivalence, $\mathcal F(f)$ is an isomorphism. The full subcategory of homotopy invariant
sheaves of spectra is denoted $\Sh_\R(\Man, \Spc)$.
\end{defn}
Constant sheaves provide good examples of homotopy invariant sheaves.

The following lemma is essentially due to Dugger~\cite{Dug01} and Morel-Voevodsky~\cite{MV99}, though they
considered space-valued sheaves. See Bunk~\cite{Bun22} for a general, model-categorical version.
\begin{lem}
The assignment $\mathcal F\mapsto\mathcal F(\pt)$ defines an equivalence $\Sh_\R(\Man, \Spc)\to\Spc$.
\end{lem}

The inclusion $\iota_\R\colon \Sh_\R(\Man, \Spc)\to\Sh(\Man, \Spc)$ admits both a left adjoint $\Lhi$ and a right
adjoint $\Rhi$. $\Rhi$ is the composition of the global sections functor $\Gamma_*\colon \Sh(\Man, \Spc)\to\Spc$
followed by the constant sheaf functor $\Gamma^*\colon\Spc\to\Sh(\Man, \Spc)$; for a formula for $\Lhi$,
see~\cite[Chapter 5]{ADH21}.
\begin{defn}
A sheaf $\mathcal F\in\Sh(\Man, \Spc)$ is \term{pure} if $\Gamma_*(\mathcal F) = 0$. The full subcategory of pure
sheaves of spectra is denoted $\Shp(\Man, \Spc)$.
\end{defn}
For example, $H\Omega^k$ is a pure sheaf. Pure sheaves tend to look like sheaves of differential forms, and contain
the ``infinitesimal'' information in a differential generalized cohomology theory.
\begin{defn}
Let $\varepsilon\colon \Rhi \Rightarrow \id$ be the counit of the adjunction $\iota_\R\dashv \Rhi$. Define a
functor $\Cyc\colon \Sh(\Man, \Spc)\to\Sh(\Man, \Spc)$ and a natural transformation $\curv\colon \id \Rightarrow
\Cyc$ by asking that $\curv\colon\id\Rightarrow \Cyc$ is the cofiber of $\varepsilon$. We call $\curv$ the
\term{curvature map} and $\Cyc(\mathcal F)$ for a sheaf $\mathcal F$ the \term{sheaf of differential cycles} of
$\mathcal F$.
\end{defn}
$\Cyc$ factors through $\Shp(\Man, \Spc)$, and is in fact left adjoint to the inclusion $\iota_{\mathrm{pu}}\colon
\Shp(\Man, \Spc)\hookrightarrow \Sh(\Man, \Spc)$.
\begin{defn}
In a similar way, let $\eta\colon \id\Rightarrow \Lhi$ be the unit of the adjunction $\Lhi\dashv \iota_\R$ and let
$\psi\colon\Def\Rightarrow \id$ be the fiber of $\eta$. Given a sheaf $\mathcal F$, $\Def(\mathcal F)$ is called
the \term{sheaf of differential deformations} of $\mathcal F$.
\end{defn}
$\Def$ is left adjoint to $\Cyc$.

The data of $\Def$, $\Cyc$, $\Lhi$, and $\Rhi$ assemble to generalizations of~\eqref{homotopy_pullback_square}
and~\eqref{the_hexagon}.
\begin{thm}[{Bunke-Nikolaus-Völkl~\cite[\S 3]{BNV16}}]
$\Sh(\Man, \Spc)$ is a \term{recollement} of its subcategories $\Sh_\R(\Man, \Spc)$ and $\Shp(\Man, \Spc)$. That
is,
\begin{enumerate}
	\item both $\iota_{\mathrm{pu}}$ and $\iota_\R$ admit left adjoints, namely $\Cyc$ and $\Lhi$;
	\item $\Cyc\circ\iota_\R\simeq 0$; and
	\item a morphism of sheaves is an equivalence if and only if both $\Cyc$ and $\Lhi$ map it to an equivalence.
\end{enumerate}
\end{thm}
\begin{cor}[Fracture square, {Bunke-Nikolaus-Völkl~\cite[Proposition 3.3]{BNV16}}]
\label{fracture_square}
There is a pullback square of natural transformations
\begin{equation}
\begin{tikzcd}
	\id & \Lhi \\
	\Cyc & \Lhi\Cyc.
	\arrow["\eta", from=1-1, to=1-2]
	\arrow["\curv"', from=1-1, to=2-1]
	\arrow["\eta", from=2-1, to=2-2]
	\arrow["\curv"', from=1-2, to=2-2]
	\arrow["\lrcorner"{anchor=center, pos=0.125}, draw=none, from=1-1, to=2-2]
\end{tikzcd}\end{equation}
\end{cor}
This is the analogue of \cref{HP_lem}: it factors an arbitrary differential generalized cohomology theory as a pullback of
something like closed forms (the pure part, in the lower left corner) and a non-differential generalized cohomology
theory (something homotopy-invariant, in the upper right corner).
\begin{cor}[Differential cohomology hexagon, {Bunke-Nikolaus-Völkl~\cite[(9)]{BNV16}}]
There is a commutative diagram of natural transformations
\begin{equation}
\begin{tikzcd}
	& \Rhi && \Lhi \\
	{\Sigma^{-1}\Lhi\Cyc} && \id && \Lhi\Cyc \\
	& \Def && \Cyc
	\arrow[from=2-1, to=1-2]
	\arrow[from=2-1, to=3-2]
	\arrow[from=1-2, to=1-4]
	\arrow[from=1-4, to=2-5]
	\arrow[from=3-2, to=3-4]
	\arrow[from=3-4, to=2-5]
	\arrow["\varepsilon"', from=1-2, to=2-3]
	\arrow["\eta"', from=2-3, to=1-4]
	\arrow["\psi", from=3-2, to=2-3]
	\arrow["\curv", from=2-3, to=3-4]
	\arrow["\lrcorner"{anchor=center, pos=0.125, rotate=45}, draw=none, from=2-3, to=2-5]
	\arrow["\lrcorner"{anchor=center, pos=0.125, rotate=45}, draw=none, from=2-1, to=2-3]
\end{tikzcd}\end{equation}
with the following properties.
\begin{enumerate}
	\item The diagonals $(\varepsilon, \curv)$ and $(\psi, \eta)$ are cofiber sequences.
	\item The top and bottom rows are once-extended cofiber sequences.
	\item Both squares are pullback squares.
\end{enumerate}
\end{cor}
Plug in a sheaf $\mathcal F$ to obtain the differential cohomology hexagon for the differential generalized
cohomology theory associated to $\mathcal F$.
\begin{rem}
This flurry of adjoints suggests that it is the presence of so many adjoints that makes the whole theory of the
differential cohomology hexagon possible. Schreiber~\cite{Sch13} takes this attitude, which he names
\term{cohesion}, and uses it to study differential cohomology in a very general setting.
\end{rem}

\subsection{Examples of differential generalized cohomology theories}
\label{gen_examples}
\begin{exm}[Ordinary differential cohomology]
For ordinary differential cohomology, apply the Eilenberg-Mac Lane functor $H$ to the Deligne complexes $\Z(n)$.
The resulting hexagon coincides with the differential cohomology hexagon from \cref{hexagon_thm}; for example,
$\Cyc(H\Z(n))\simeq H\R(n)$, recovering the sheaf of closed forms.
\end{exm}
\begin{exm}[Differential $K$-theory]
\label{diffK}
Differential $K$-theory was first studied by Freed~\cite{Fre00} and Freed-Hopkins~\cite{FH00} for applications in
string theory, with related objects considered earlier by Gillet-Soulé~\cite{GS90} and Lott~\cite{Lot00}.
Hopkins-Singer~\cite[\S 4.4]{HS05} gave the first comprehensive construction of differential $K$-theory, and
additional constructions have been given by
Klonoff \cite{Klo08},
Bunke-Schick \cite[\S 2]{BS09},
Simons-Sullivan \cite{SS10},
Bunke-Nikolaus-Völkl \cite[\S 6]{BNV16},
Schlegel \cite[\S 4.2]{Sch13a},
Tradler-Wilson-Zenalian \cite{TWZ13, TWZ16},
Hekmati-Murray-Schlegel-Vozzo \cite{HMSV15},
Park \cite{Par17},
Gorokhovski--Lott \cite{GL18},
Schlarmann \cite{Sch19},
Cushman~\cite{Cus21},
Park-Parzygnat-Redden-Stoffel \cite{PPRS22},
Gomi-Yamashita~\cite{GY23}, and
Lee-Park~\cite{LP23}.
See Bunke--Schick \cite{BunksSchick} for a survey.

The idea of differential $K$-theory is to use the Chern character as the source of differential form information
refining a $K$-theory class. Let $A\coloneqq \KU^*(\pt)\cong\Z[t,t^{-1}]$, with $\abs t = 2$. The \term{Chern
character} is the map of spectra sending $\KU$ to its tensor product with $\R$:
\begin{equation}
	\ch\colon \KU\longrightarrow \KU\wedge H\R \simeq H(\R\otimes A).
\end{equation}
Fix $n\in\Z$, though only the value of $n\bmod 2$ will matter in the end, due to Bott periodicity. Then define a
$K$-theoretic analogue of the Deligne complex $\KU(n)$ as the homotopy pullback
\begin{equation}\label{HS_K_diag}
\begin{tikzcd}
	{\KU(n)} & \KU \\
	{H\tau_{\ge 0}(\R(n)\otimes A)} & {H(\R\otimes  A)}.
	\arrow["t", from=2-1, to=2-2]
	\arrow["\ch"', from=1-2, to=2-2]
	\arrow[from=1-1, to=1-2]
	\arrow[from=1-1, to=2-1]
	\arrow["\lrcorner"{anchor=center, pos=0.125}, draw=none, from=1-1, to=2-2]
\end{tikzcd}\end{equation}
In the lower left corner, $\tau_{\ge 0}$ means taking the nonnegatively graded parts only (since $A$ can contribute
negative grading). This sheaf consists of closed $\R\otimes A$-valued forms whose degrees, possibly shifted by
multiplication by a power of $t$, are nonnegative and of the same parity as $n$. The reason for this complicated
object is that the Chern character associated to a connection on a vector bundle is a form of this type.

The \term{differential $K$-theory} groups $\check K^n(M)$ are the hypercohomology groups $H^n(M; \KU(n))$. They are
$2$-periodic, like for ordinary $K$-theory. $\check K^0(M)$ is naturally isomorphic to the group completion of the
commutative monoid of vector bundles with connection on $M$.

We can then fill in the rest of the hexagon for differential $K$-theory. This diagram was first constructed by
Simons-Sullivan~\cite{SS10}:
\begin{equation}\label{kth_hexagon}
\begin{tikzcd}
	& {K^n(M;\R/\Z)} && {K^n(M)} \\
	{H^{n-1}(M;\R\otimes A)} && {\check K^n(M)} && {H^n(M; \R\otimes A)} \\
	& {\Omega^{n-1}(M;\R\otimes A)/\Im(\d)} && {\Omega^n(M; \R\otimes A)_{\mathit{c\ell}}}
	\arrow["\cc", from=2-3, to=1-4]
	\arrow["{\ch^\nabla}"', from=2-3, to=3-4]
	\arrow["\ch", from=1-4, to=2-5]
	\arrow[from=3-4, to=2-5]
	\arrow[from=2-1, to=1-2]
	\arrow[from=1-2, to=1-4]
	\arrow[from=1-2, to=2-3]
	\arrow[from=2-1, to=3-2]
	\arrow["\d", from=3-2, to=3-4]
	\arrow[from=3-2, to=2-3]
	\arrow["\lrcorner"{anchor=center, pos=0.125, rotate=45}, draw=none, from=2-1, to=2-3]
	\arrow["\lrcorner"{anchor=center, pos=0.125, rotate=45}, draw=none, from=2-3, to=2-5]
\end{tikzcd}\end{equation}
The story is still roughly similar to the hexagon for ordinary differential cohomology, but there is some new
notation.
\begin{itemize}
	\item $K^*(\text{--};\R/\Z)$ is \term{$K$-theory with $\R/\Z$ coefficients}, the generalized cohomology theory
	represented by the spectrum which is the cofiber of $\ch\colon\KU\to \KU\wedge H\R\simeq H(\R\otimes A)$. This
	theory first appears in Atiyah-Patodi-Singer~\cite[\S 5]{APS3}, who attribute it to Segal. The
	long exact sequence in cohomology induced by the cofiber sequence $\KU\to\KU\wedge H\R\to
	\KU(\text{--};\R/\Z)$, which is a $K$-theoretic analogue of the $\Z\to\R\to\R/\Z$ Bockstein long exact
	sequence, is the upper long exact sequence in~\eqref{kth_hexagon}.
	\item $\cc$ is the characteristic class map, which is the topmost map in~\eqref{HS_K_diag}.
	\item $\ch^\nabla$ is the version of the Chern character which takes in a vector bundle with connection and
	produces a closed form. This is the curvature map for differential $K$-theory.
\end{itemize}
\end{exm}
\begin{exm}[Differential $\KO$-theory]
Like differential $K$-theory, differential $\KO$-theory was first studied by Freed~\cite{Fre00} and
Freed-Hopkins~\cite{FH00}; Grady-Sati~\cite{GS19, GS21} were the first to comprehensively study differential
$\KO$-theory, and Cushman~\cite{Cus21} and Gomi-Yamashita~\cite{GY23} provide additional constructions.

The real version of \cref{diffK} is completely analogous. Instead of using the Chern character, one uses its real
analogue (sometimes called the \term{Pontrjagin character})
\begin{equation}
	\ph\colon\KO\longrightarrow \KO\wedge H\R\simeq H(\R\otimes B),
\end{equation}
where $B\cong\Z[t,t^{-1}]$ with $\abs t = 4$.\footnote{$B$ is not isomorphic to $\KO^*(\pt)$. When we tensor with
$\R$, this discrepancy goes away.} One succinct way to define both $\ph$ and its form-level version $\ph^\nabla$
for a real vector bundle with connection is: first complexify, then take the Chern character. The result a priori
lands in $\R\otimes A$-valued forms (resp.\ cohomology), but in fact factors through $\R\otimes B$-valued forms
(resp.\ cohomology).

Thus we have Deligne-type complexes, now depending on $n\bmod 8$:
\begin{equation}\label{HS_KO_diag}
\begin{tikzcd}
	{\KO(n)} & \KO \\
	{H\tau_{\ge 0}(\R(n)\otimes B)} & {H(\R\otimes  B)},
	\arrow["t", from=2-1, to=2-2]
	\arrow["\ph"', from=1-2, to=2-2]
	\arrow[from=1-1, to=1-2]
	\arrow[from=1-1, to=2-1]
	\arrow["\lrcorner"{anchor=center, pos=0.125}, draw=none, from=1-1, to=2-2]
\end{tikzcd}\end{equation}
and the differential $\KO$-cohomology hexagon:
\begin{equation}
\begin{tikzcd}
	& {\KO^n(M;\R/\Z)} && {\KO^n(M)} & {} \\
	{H^{n-1}(M;\R\otimes B)} && {\check \KO^n(M)} && {H^n(M; \R\otimes B),} \\
	& {\Omega^{n-1}(M;\R\otimes B)/\Im(\d)} && {\Omega^n(M; \R\otimes B)_{\mathit{c\ell}}}
	\arrow["\cc", from=2-3, to=1-4]
	\arrow["{\ph^\nabla}"', from=2-3, to=3-4]
	\arrow["\ph", from=1-4, to=2-5]
	\arrow[from=3-4, to=2-5]
	\arrow[from=2-1, to=1-2]
	\arrow[from=1-2, to=1-4]
	\arrow[from=1-2, to=2-3]
	\arrow[from=2-1, to=3-2]
	\arrow[from=3-2, to=3-4]
	\arrow[from=3-2, to=2-3]
	\arrow["\lrcorner"{anchor=center, pos=0.125, rotate=45}, draw=none, from=2-3, to=2-5]
	\arrow["\lrcorner"{anchor=center, pos=0.125, rotate=45}, draw=none, from=2-1, to=2-3]
\end{tikzcd}\end{equation}
where as usual $\check{\KO}^n(M)$ is the $n^{\mathrm{th}}$ cohomology of $M$ valued in $\KO(n)$.
\end{exm}
\begin{rem}[Some more examples]
Though differential $K$- and $\KO$-theory are the most commonly studied differential generalized cohomology
theories, several others appear in the literature.
\begin{enumerate}
	\item \term{Supercohomology} $\mathit{SH}$, defined by Freed~\cite[\S 1]{Fre08} and Gu-Wen~\cite{GW14}, is the
	spectrum with $\pi_0(\mathit{SH})\cong\Z$, $\pi_2(\mathit{SH})\cong\Z/2$, and the unique nontrivial Postnikov
	invariant connecting them.\footnote{Sometimes $\mathit{SH}$ is called \term{restricted supercohomology} to
	contrast with \term{extended supercohomology}, a different spectrum studied by Kapustin-Thorngren~\cite{KT17}
	and Wang-Gu~\cite{WG20}. See~\cite[\S 5.3, 5.4]{GJF19}.} Freed-Neitzke~\cite{FN22, FN23} introduce a
	differential refinement of this theory for the purpose of studying classical spin Chern-Simons theory.
	\item Differential refinements of algebraic $K$-theory spectra appear in work of Bunke-Tamme~\cite{BT15, BT16},
	Bunke~\cite{Bun18b, Bun18a}, Schrade~\cite{Sch18}, Bunke-Gepner~\cite{BG21}, and
	Park-Parzygnat-Redden-Stoffel~\cite{PPRS22}, where among other things they are applied to construct a
	topological version of Beĭlinson's regulator homomorphisms.
	\item A complex-analytic differential refinement of $\mathit{MU}$, the spectrum representing complex cobordism,
	appears in work of Quick and collaborators~\cite{HQ15, Qui16, Qui19, HQ23a, HQ23b, KQ23}; another differential
	cobordism theory appears in work of Bunke-Schick-Schröder-Weithaup~\cite{BSSW09}. See also
	Grady-Sati~\cite{GS17} for a closely related construction.
\end{enumerate}
\end{rem}

\section{Applications in physics}
	\label{s_physics}
	Closed differential forms are commonplace in the classical theory of electromagnetism, encoding quantities such as
the field strength. Passing to the quantum theory amounts to choosing integrality data for the de Rham classes of
these forms --- in other words, lifting them to differential cocycles. We will discuss this story in this section,
where it also leads to the original motivation for differential generalized cohomology theories
(\S\ref{diffgenph}).

\subsection{Dirac quantization in electromagnetism}
For the first part of this section, we follow Freed~\cite{Fre00}; see also~\cite[Chapter 21]{ADH21}.

Let us go over the basic objects of nonrelativistic classical electromagnetic theory in three-dimensional space,
which for us will be an oriented Riemannian $3$-manifold $Y$ with empty boundary. Let $X\coloneqq \R\times Y$, with
the Lorentzian metric $\d t^2 - g$, where $t$ is the $\R$-coordinate.

One may be used to thinking of an electric field as a vector field, representing at each point the magnitude and
direction of force exerted on a unit test charge. We will use the metric to pass between $TY$ and $T^*Y$ and
describe the electric field as a $1$-form $E\in\Omega^1(Y)$. For the magnetic field, it is helpful to instead pass
through the Hodge star and obtain a $2$-form $B\in\Omega^2(Y)$. The \term{charge density} $\rho_c$ is a compactly
supported differential $3$-form, and the \term{(electric) current} $J_E$ is a compactly supported $2$-form.

The \term{field strength} is $F\coloneqq B - \d t\wedge E\in\Omega^2(X)$, and let $j_E\coloneqq \rho_E - \ud
t\wedge J_E\in\Omega^3_c(X)$. Maxwell's equations can be concisely expressed in terms of $F$ and $j_E$:
\begin{equation}
\begin{aligned}
	\d F &= 0\\
	\d {\star F} &= j_E.
\end{aligned}
\end{equation}
If there is a \term{magnetic current} $j_B\in\Omega^3(X)$, we modify the first equation to $\d F = j_B$.

One can then use these forms to write down a Lagrangian action, compute quantities such as the total charge, and so
on. The total charge $\overline Q$ is a cohomological object, in fact --- it is the de Rham class of $j_E$ in
$H^3_c(Y;\R)$. There is an analogous total magnetic charge.

Quantization tells us that the total charge ought to be discrete --- for example, if $Y = \R^3$,
$H_c^3(\R^3;\R)\cong\R$, and we assume the total electric charge is some integer multiple of a unit charge $q_E$.
In general, we postulate that the charge must be in the image of the map $H_c^3(Y; q_E\Z)\to H_c^3(Y;\R)$, and
likewise for a unit magnetic charge $q_B$.

So the electric charge is a closed form with what looks like data of a lift of its de Rham class to $q_E\Z$-cohomology. This suggests:
\begin{ansatz}
\label{first_ansatz}
Objects represented by closed differential forms in a classical theory of physics should be represented by cocycles
for ordinary differential cohomology in the corresponding quantum theory.
\end{ansatz}
We use cocycles, rather than cohomology classes, in order to obtain something which sheafifies, part of the
principle of locality of quantum field theory.

In general, differential forms represent plenty of objects in field theories. Notably, they are gauge fields for
abelian gauge groups, including for ``higher gauge theory'' where the gauge group is a categorification of the
circle group and one uses (higher) gerbes instead of principal bundles with connection.

\subsection{Quantizing in more general cohomology theories}
\label{diffgenph}

String theory teaches us a striking lesson: that for some differential forms, the natural home for the fields in
the quantized theory is a differential generalized cohomology theory. Typically this is differential $K$- or
$\KO$-theory, but choosing the correct theory is more of an art than a science and there are different proposals
using different differential generalized cohomology theories.

For example, consider type IIB string theory on a $10$-manifold $X$. There is a $3$-form field $B$, which as above
should be upgraded to a cocycle for $\check B\in\check H^3(X;\Z)$. For now, assume this field is zero;\footnote{If
this field is nonzero, one should repeat this discussion with \term{twisted differential $K$-theory}.} then there
are several forms called \term{Ramond-Ramond field strengths} $G_i\in\Omega^i(X)$, where $i = 1,3,5,7,9$. These
field strengths satisfy related integrality conditions implying that they are the Chern character of a cocycle for
$\check K^1(X)$, so we postulate that the Ramond-Ramond field \emph{is} a cocycle for differential $K$-theory. See
Freed~\cite[\S 3]{Fre02} for further discussion. Other examples include $\check K^0$ appearing in type IIA string
theory, $\check{\KO}^*$ in type I string theory, and the type II $B$-field lifting to a differential refinement of
a Postnikov truncation of $\mathrm{Pic}(\KU)$-cohomology, as described by Distler-Freed-Moore~\cite{DFM11a, DFM11b}.

See also \cite{BM06a, BM06b, DFM07, FMS07a, FMS07b, Fre08, Sat10, SZ10, Sat11, SSS12, KM13, KV14, DMDR14, FSS15,
FR16, GS19, Sat19, FRRB20, SS23, SS23a} for more examples of quantization in differential generalized cohomology
theories. Of particular note is ``hypothesis H'' of Fiorenza--Sati--Schreiber \cite{Sat18, FSS19, FSS20a} proposing
that the $C$-field in M-theory is quantized using twisted differential cohomotopy; work of Fiorenza, Sati,
Schreiber, and their collaborators \cite{FSS19, SS19, FSS20a, FSS20b, FSS20c, GS20, SS20b, SS20a, SS20c, BSS21,
SS21, SS23a}, as well as Roberts \cite{Rob20}, explores this hypothesis and its consequences.

%
%
%

\section{Further reading}
	\label{further_reading}
The book~\cite{ADH21} is an introduction to differential cohomology with much the same attitude as the current
article; we also recommend the other book-length introductions~\cite{BB14, Bunke:notes,
Sch13}. \cite{HS05}, is a research article that we also recommend as a
book-length introduction.

One application of differential generalized cohomology in physics that we did not get into is the classification of
reflection-positive invertible field theories, conjectured by Freed-Hopkins~\cite{FH21} and proven by
Grady~\cite{Gra23}. See Freed-Hopkins~\cite{FH21} and Freed~\cite{Fre19} for more on this conjecture,
and~\cite[Chapter 22]{ADH21} for a review, and see Davighi-Gripaios-Randal-Williams~\cite{DGR20}, Yamashita-Yonekura~\cite{YY23}, and Yamashita~\cite{Yam23a,
Yam23b} for some related work.







\end{document}